\documentclass[envcountsame]{llncs}
\usepackage{gastex}
\usepackage{amsmath}
\usepackage{amssymb}

\def\blfootnote{\xdef\@thefnmark{}\@footnotetext}
\long\def\symbolfootnote[#1]#2{\begingroup%
\def\thefootnote{\fnsymbol{footnote}}\footnote[#1]{#2}\endgroup}

\newenvironment{Theorem}{\begin{theorem}}{\end{theorem}}
\newenvironment{Proposition}{\begin{proposition}}{\end{proposition}}
\newenvironment{Lemma}{\begin{lemma}}{\end{lemma}}
\newenvironment{Corollary}{\begin{corollary}}{\end{corollary}}
\newenvironment{Remark}{\begin{remark}\rm}{\qed\end{remark}}

\newenvironment{Proof}{\begin{proof}}{\qed\end{proof}}

\newenvironment{Example}{\begin{example}\rm}{\qed\end{example}}

\newcommand{\xor}{\oplus}
\newcommand{\pmatch}{\mathrm{pm}}
\newcommand{\tpart}{\mathrm{pair}}
\newcommand{\sgn}{\mathrm{sgn}}
\newcommand{\two}{GF(2)}
\newcommand{\pair}[2]{{#1}{#2}}
\newcommand{\sub}[1]{\langle #1\rangle}
\newcommand{\HJ}[1]{}
\newcommand{\doi}[1]{doi: http://dx.doi.org/#1}

\begin{document}

\title{Pivots, Determinants, and Perfect Matchings of Graphs}
\author{Robert Brijder\inst{1}%
        \thanks{corresponding author: \email{rbrijder@liacs.nl}}
        \and
        Tero Harju\inst{2} \and
        Hendrik Jan Hoogeboom\inst{1}}

\institute{Leiden Institute of Advanced Computer Science, Universiteit Leiden,\\
Niels Bohrweg 1, 2333 CA Leiden, The Netherlands
\and
Department of Mathematics,
University of Turku, FI-20014 Turku, Finland
}

\maketitle
\begin{abstract}
We give a characterization of the effect of sequences of pivot
operations on a graph by relating it to determinants of adjacency
matrices. This allows us to deduce that two sequences of pivot
operations are equivalent iff they contain the same set $S$ of
vertices (modulo two). Moreover, given a set of vertices $S$, we
characterize whether or not such a sequence using precisely the
vertices of $S$ exists. We also relate pivots to perfect matchings
to obtain a graph-theoretical characterization. Finally, we consider
graphs with self-loops to carry over the results to sequences
containing both pivots and local complementation operations.
\end{abstract}

\section{Introduction}
The operation of local complementation in an undirected graph takes
the neighbourhood of a vertex in the graph and replaces that
neighbourhood by its graph complement. The related operation of edge
local complementation, here called pivoting, can be defined in terms
of local complementation. It starts with an edge in the graph and
toggles edges based on the way its endpoints are connected to the
endpoints of the pivot-edge.

The operations are connected in a natural way to overlap graphs
(also called circle graphs \cite{gavril}). Given a finite set of
chords of a circle, the overlap graph contains a vertex for each
chord, and two vertices are connected if the corresponding chords
cross. Taking out a piece of the perimeter of the circle delimited
by the two endpoints of a chord, and reinserting it in reverse,
changes the way the cords intersect, and hence changes the
associated overlap graph. The effect of this reversal on the overlap
graph can be obtained by a local complementation on the vertex
corresponding to the chord. Similarly, interchanging two pieces of
the perimeter of the circle, each starting at the different
endpoints of one common chord and ending at the endpoints of
another, can be modelled by a pivot in the overlap graph.

Overlap graphs naturally occur in theories of genetic rearrangements
\cite{pevzner,ciliates}, but local complementation and edge local
complementation operations are applied in many settings, like the
relationships between Eulerian tours, equivalence of certain
codes\cite{danielsen}, rank-width of graphs\cite{oum}, and quantum
graph states\cite{nest}.

In the present paper we are interested in sequences of pivots in
arbitrary simple graphs. In defining a single pivot one usually
distinguishes three disjoint neighbourhoods in the graph, and edges
are updated according to the neighbourhoods to which the endpoints
belong. Describing the effect of a sequence of pivot operations in
terms of neighbourhood connections is involved -- the number of
neighbourhoods to consider grows exponentially in the size of the
sequence.

It turns out that by considering determinants of adjacency matrices
(in the spirit of \cite{geelen}) we can effectively describe the
effect of sequences of pivot operations. Subsequently, we relate
this to perfect matchings, a perfect matching is a set of edges that
forms a partition of the set of vertices, to obtain a
graph-theoretical characterization. A direct proof of the
characterization in terms of perfect matchings is given in the
appendix. We obtain the surprising result that the connection
between two vertices after a series of pivots directly depends on
the number (modulo two) of perfect matchings in the subgraph induced
by the two vertices and the vertices of the pivot-edges (with
`multiplicity' if vertices occur more than once).

As an immediate consequence we obtain that the result of a sequence
of pivots, provided all pivot operations are defined, i.e., based on
an edge in the graph to which they are applied, does not depend on
the order of the pivots, but only on the nodes involved (plus their
cardinality modulo 2). Also, we show that for any applicable
sequence of pivot there exists an equivalent \emph{reduced} sequence
where each node appears at most once in the sequence. Finally, we
consider the case where graphs can have self-loops, and generalize
the results for sequences of pivots to sequences having both local
complementation operations and pivots.

\section{Preliminaries} \label{sec:preliminaries}

Usually we write $\pair xy$ for the pair $\{x,y\}$.

We use $\xor$ to denote both the logical exclusive-or as well as the
related operation of symmetric set difference. The operation is
$\xor$ associative: the exclusive or over a sequence of Booleans is
true iff an odd number of the arguments is true.

Let $A$ be a $V\times V$ matrix. For a set $X\subseteq V$ we use
$A\sub{X}$ to denote the submatrix induced by $X$, which keeps the
rows and columns indexed by $X$.

The determinant of $A$ is defined as
$\det (A) = \sum_{\sigma\in \Pi(V)} \sgn(\sigma)
\prod_{u\in V} a_{u,\sigma(u)}$,
where $\Pi(V)$ is the set of permutations of $V$, and $\sgn(\sigma)$
is the sign (or parity) of the permutation, which is well defined
after choosing an ordering on $V$. We will mainly consider the
determinant over $\two$, i.e., modulo 2, where the signs do not
matter. The determinant of the empty matrix is considered to be $1$
(contributed by the empty permutation).

\paragraph{Graphs.}
The graphs we consider here are simple (undirected and without loops
and parallel edges).
For graph $G=(V,E)$ we use $V(G)$ and $E(G)$ to denote its set of
vertices $V$ and set of edges $E$, respectively.

We define $x \sim_G y$ if either $\pair xy \in E$ or $x=y$.
For $X \subseteq V$, we denote the subgraph of $G$
induced by $X$ as $G\sub{X}$.
Let $N_G(v) = \{w \in V \mid \pair vw \in E\}$
denote the neighbourhood of vertex $v$ in graph $G$.

With a graph $G$ one associates its adjacency matrix $A(G)$, which
is a $V\times V$ $(0,1)$-matrix $\left(a_{u,v}\right)$ with $a_{u,v}
= 1$ iff $\pair uv \in E$. Obviously, for $X\subseteq V$,
$A(G\sub{X}) = A(G)\sub{X}$.

By the determinant of graph $G$, denoted $\det G$,
we will mean the determinant $\det A(G)$ of its adjacency
matrix, computed over $\two$.

\section{Pivot Operation} \label{sec_pivots}

Let $G = (V,E)$ be a graph. The graph obtained by \emph{local
complementation} at $u \in V$ on $G$, denoted by $G * u$, is the
graph that is obtained from $G$ by complementing the edges in the
neighbourhood $N_G(u)$ of $u$ in $G$. Using a logical expression we
can write for $G*u $ the definition $\pair xy \in E(G*u)$ iff
$(xy\in E)\xor (\pair xu\in E \land yu\in E)$.

For a vertex $x$ consider its closed neighbourhood
$N'_G(x)= N_G(x)\cup \{x\} = \{ y\in V_G \mid x\sim_G y \}$.
The edge $uv$ partitions the vertices of $G$
connected to $u$ or $v$ into three sets
$V_1 = N'_G(u) \setminus N'_G(v)$,
$V_2 = N'_G(v) \setminus N'_G(u)$,
$V_3 = N'_G(u) \cap N'_G(v)$.
Note that $u,v \in V_3$.

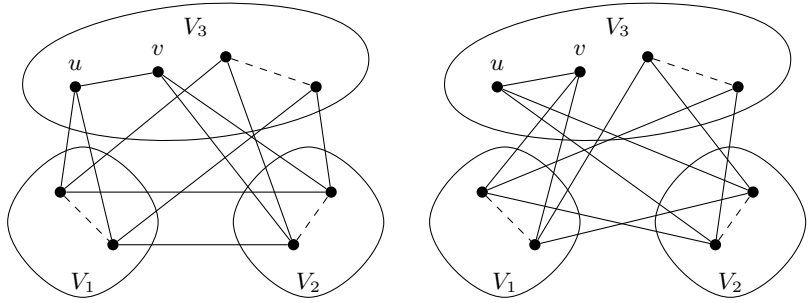
\begin{figure}[t]
\centerline{\unitlength 1.0mm
\begin{picture}(55,42)(0,1)
\drawccurve(02,28)(25,21)(48,32)(25,39)
\drawccurve(00,10)(10,00)(20,10)(10,20)
\drawccurve(30,10)(40,00)(50,10)(40,20)
\gasset{AHnb=0,Nw=1.5,Nh=1.5,Nframe=n,Nfill=y}
\gasset{ExtNL=y,NLdist=1.5,NLangle=90}
\put(10,02){\makebox(0,0)[cc]{$V_1$}}
\put(40,02){\makebox(0,0)[cc]{$V_2$}}
\put(25,36){\makebox(0,0)[cc]{$V_3$}}
  \node(u)(09,28){$u$}
  \node(v)(20,30){$v$}
  \node(uu)(29,32){}
  \node(vv)(41,28){}
  \node(u1)(7,14){}
  \node(u2)(14,7){}
  \node(v1)(38,7){}
  \node(v2)(43,14){}
  \drawedge(u,v){}
  \drawedge(u,u1){}
  \drawedge(u,u2){}
  \drawedge(v,v1){}
  \drawedge(v,v2){}
  \drawedge(u1,v2){}
  \drawedge(u2,v1){}
  \drawedge[dash={1}0](v1,v2){}
  \drawedge[dash={1}0](u1,u2){}
  \drawedge[dash={1}0](uu,vv){}
  \drawedge(uu,u1){}
  \drawedge(vv,u2){}
  \drawedge(uu,v1){}
  \drawedge(vv,v2){}
\end{picture}
\begin{picture}(55,42)(0,1)
\drawccurve(02,28)(25,21)(48,32)(25,39)
\drawccurve(00,10)(10,00)(20,10)(10,20)
\drawccurve(30,10)(40,00)(50,10)(40,20)
\gasset{AHnb=0,Nw=1.5,Nh=1.5,Nframe=n,Nfill=y}
\gasset{ExtNL=y,NLdist=1.5,NLangle=90}
\put(10,02){\makebox(0,0)[cc]{$V_1$}}
\put(40,02){\makebox(0,0)[cc]{$V_2$}}
\put(25,36){\makebox(0,0)[cc]{$V_3$}}
  \node(u)(09,28){$u$}
  \node(v)(20,30){$v$}
  \node(uu)(29,32){}
  \node(vv)(41,28){}
  \node(u1)(7,14){}
  \node(u2)(14,7){}
  \node(v1)(38,7){}
  \node(v2)(43,14){}
  \drawedge(u,v){}
  \drawedge(v,u1){}
  \drawedge(v,u2){}
  \drawedge(u,v1){}
  \drawedge(u,v2){}
  \drawedge(u1,v1){}
  \drawedge(u2,v2){}
  \drawedge[dash={1}0](v1,v2){}
  \drawedge[dash={1}0](u1,u2){}
  \drawedge[dash={1}0](uu,vv){}
  \drawedge(uu,u2){}
  \drawedge(vv,u1){}
  \drawedge(uu,v2){}
  \drawedge(vv,v1){}
\end{picture}%
}
\caption{Pivoting $\pair uv$.
 Connection $\pair xy$ is toggled if $x\in V_i$ and $y\in V_j$ with
 $i\neq j$. Note $u$ and $v$ are connected to all vertices in $V_2$,
 these edges are omitted in the diagram. The operation does not
 effect edges adjacent to vertices outside the sets
 $V_1,V_2,V_3$.}%
\label{fig:pivot}
\end{figure}

Let $\pair uv \in E(G)$. The graph obtained from $G$ by
\emph{pivoting} $\pair uv$, denoted by $G[uv]$, is constructed by
`toggling' all edges between different $V_i$ and $V_j$: for $\pair
xy$ with $x\in V_i$ and $y\in V_j$ ($i\neq j$): $\pair xy \in E(G)$
iff $\pair xy \notin E(G[uv])$, see Figure~\ref{fig:pivot}. The
remaining edges remain unchanged. \footnote{In defining this
operation usually the
 description adds the rule that the vertices $u$ and $v$
 are swapped. Here this is avoided by including $u$ and $v$
 in the set $V_3$.}

It turns out that $G[uv]$ equals $G *u*v*u = G *v*u*v$.

\begin{Example}\label{ex:overlap}
We start with six segments, of which the relative positions of
endpoints can be
represented by the
string
$3\; 5\; 2\; 6\; 5\; 4\; 1\; 3\; 6\; 1\; 2\; 4 $.

The `entanglement' of these intervals can be represented
by the overlap graph to the left
in Figure~\ref{fig:overlap}.
When we pivot on the edge $\pair 23$ we obtain the graph
to the right.

This new graph is the overlap graph of
$\underline{3\; 6\; 1\; 2}\; 6\; 5\; 4\; 1\; \underline{3\; 5\; 2}\; 4 $.

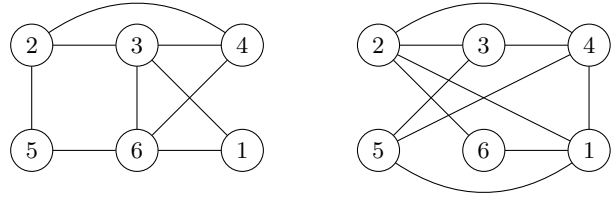
\begin{figure}
\centerline{\unitlength 0.7mm%
\begin{picture}(50,40)
\gasset{AHnb=0,Nw=1.5,Nh=1.5,Nframe=n,Nfill=y}
\gasset{AHnb=0,Nw=8,Nh=8,Nframe=y,Nfill=n}
  \node(2)(05,25){2}
  \node(3)(25,25){3}
  \node(4)(45,25){4}
  \node(5)(05,05){5}
  \node(6)(25,05){6}
  \node(7)(45,05){1}
  \drawedge(2,3){}
  \drawedge[curvedepth=8](2,4){}
  \drawedge(2,5){}
  \drawedge(3,4){}
  \drawedge(3,6){}
  \drawedge(3,7){}
  \drawedge(4,6){}
  \drawedge(5,6){}
  \drawedge(6,7){}
\end{picture}\hspace{1cm}
\begin{picture}(50,40)
\gasset{AHnb=0,Nw=1.5,Nh=1.5,Nframe=n,Nfill=y}
\gasset{AHnb=0,Nw=8,Nh=8,Nframe=y,Nfill=n}
  \node(2)(05,25){2}
  \node(3)(25,25){3}
  \node(4)(45,25){4}
  \node(5)(05,05){5}
  \node(6)(25,05){6}
  \node(7)(45,05){1}
  \drawedge(2,3){}
  \drawedge[curvedepth=8](2,4){}
  \drawedge(2,6){}
  \drawedge(2,7){}
  \drawedge(3,4){}
  \drawedge(3,5){}
  \drawedge(4,5){}
  \drawedge[curvedepth=-8](5,7){}
  \drawedge(6,7){}
  \drawedge(4,7){}
\end{picture}%
} \caption{A graph $G$ and its pivot $G[\pair 23]$, cf.
Example~\ref{ex:overlap}.}\label{fig:overlap}
\end{figure}
\end{Example}

In order to derive properties of pivoting in an algebraic way,
rather than using combinatorial methods in graphs,
Oum \cite{oum} shows that $G[uv]$ can be described
using a logical formula.
It turns out that the expression can be stated
elegantly in terms of $\sim_G$ rather than in terms of $E(G)$.

\begin{Lemma}\label{lem:oum}
Let $G$ be a graph, and let $\pair uv\in E(G)$.
Then $G[uv]$ is defined by the expression
$$
x \sim_{G[uv]} y = x \sim_G y \xor
((x \sim_G u) \wedge (y \sim_G v)) \xor
((x \sim_G v) \wedge (y \sim_G u)).
$$
for all $x,y \in V(G)$.\qed
\end{Lemma}

\paragraph{Pivots and matrices.}
In a 1997 paper \cite{geelen} on unimodular $(0,1)$-matrices, Geelen
defines a general pivot operation on matrices that is defined for
subsets of the indices (thus not only for edges) which turns out to
extend the classic pivot operation introduced above.

Let $A$ be a $V$ by $V$ $(0,1)$-matrix, and let $X\subseteq V$ be
such that $\det A\sub{X} \neq 0$, then the \emph{pivot} of $A$ by
$X$, denoted by $A*X$\footnote{The local complementation operation
G*u differs from $A*\{u\}$ where $A$ is the adjacency matrix, see
Section~\ref{sec_self_loops}}, is defined as follows. If $P =
A\sub{X}$ and $A = \left(
\begin{array}{c|c}
P & Q \\
\hline R & S
\end{array}
\right)$, then
$$
A*X = \left(
\begin{array}{c|c}
-P^{-1} & P^{-1} Q \\
\hline R P^{-1} & S - R P^{-1} Q
\end{array}
\right).
$$

Based on a similar operation from \cite{tucker} (see also
\cite[p.230]{cottle}), the following basic result can be obtained,
see \cite[Theorem~2.1]{geelen} and \cite{geelenphd} for a full
proof.

\begin{Proposition}\label{prop:geelen}
Let $A$ be a $V\times V$ matrix, and let
$X\subseteq V$ be such that $\det A\sub{X} \neq 0$.
Then, for $Y \subseteq V$,
\[
    \det (A*X)\sub{Y} = \pm \det A\sub{X \xor Y} / \det A\sub{X}
\]\qed
\end{Proposition}

We will apply this result to our (edge) pivots in graphs.
Let $A$ be the adjacency matrix of graph $G$.
we start by observing that for vertices $u\neq v$,
$\pair uv$ is an edge in $G$ iff
the submatrix $A\sub{\pair uv}$ is of the form
$\left(
\begin{array}{cc}
0 & 1 \\
1 & 0
\end{array}
\right)
$
or equivalently
$\det G\sub{\pair uv} = 1$.

If $\pair uv$ is an edge in $G$, then
(after rearranging rows and columns)
$A$ can be written in the form
\[
A  = \left(
\begin{array}{c|c|c}
0 & 1 & \chi_u^T \\\hline
1 & 0 & \chi_v^T \\\hline
\chi_u & \chi_v & A\sub{V-u-v}
\end{array}
\right)
\]
where $\chi_u$ is the column vector belonging to $u$ without
elements $a_{uu}$ and $a_{vu}$, and, for vector $x$, $x^T$ is the
transpose of $x$.

As  $\det A\sub{\pair uv} \neq 0$,
the operation pivot
$A * \pair uv$ of \cite{geelen} is
well defined. It equals the following matrix
which in fact is the matrix of $G[\pair uv]$:
the component $(\chi_v\chi_u^T + \chi_u\chi_v^T)$
in the matrix has the same functionality as
the expression $((x \sim_G u) \wedge (y \sim_G v)) \xor
((x \sim_G v) \wedge (y \sim_G u))$
from the characterization of Oum, Lemma~\ref{lem:oum}.
\[
A * \pair uv  = \left(
\begin{array}{c|c|c}
0 & 1 & \chi_v^T \\\hline
1 & 0 & \chi_u^T \\\hline
\chi_v & \chi_u & A\sub{V-u-v} - (\chi_v\chi_u^T + \chi_u\chi_v^T)
\end{array}
\right)
\]

We now rephrase the result cited from \cite{geelen},
Proposition~\ref{prop:geelen} above,
for pivots in graphs (and where the computations
are over $\two$).
It will be the main tool in our paper.


\begin{Theorem}\label{thm:geelen}
Let $G$ be a graph, and let $\pair uv \in E(G)$. Then, for $Y
\subseteq V(G)$,
\[
    \det ((G[\pair uv])\sub{Y}) = \det (G\sub{Y \xor \{u,v\}})
\]\qed
\end{Theorem}

It is noted in Little~\cite{little} that over $\two$ the $\det (G) =
0$ has a graph interpretation: $\det (G) = 0$ iff there exists a
non-empty set $S \subseteq V(G)$ such that every $v \in V(G)$ is
adjacent to an even number of vertices in $S$. Indeed, $S$
represents a linear dependent set of rows modulo $2$.

Finally note that $x \sim_G y$ iff $\det (G[\{x\}\xor\{y\}]) = 1$.
Indeed, if $x=y$, then $\det (G[\emptyset]) = 1$, and if $x \not=
y$, then $xy$ is an edge iff $\det (G[xy]) = 1$.

\section{Sequences of Pivots} \label{sec_seq_pivots}

In this section we study series of pivots that are applied
consecutively to a graph.
It is shown that using determinants there is an
elegant formula that describes whether a certain pair of
vertices is adjacent in the final resulting graph.
From this result we then conclude that the effect of
a sequence of pivots only depends on the vertices
involved, and not on the order of the operations.
Without determinants, using combinatorical argumentations
on graphs, this result seems hard to obtain.

A sequence of pivoting operations $\varphi
= [v_1v_2] [v_3v_4] \cdots [v_{n-1}v_n]$
is \emph{applicable} if each pair $[v_{i}v_{i+1}]$ in the sequence
corresponds to an edge $\pair {v_{i}}{v_{i+1}}$
in the graph obtained at the time of
application.
For such a sequence we
define $\sup(\varphi) = \bigoplus_i \{v_i\}$, the set of vertices that
occur an odd number of times in the sequence of operations.
This is called the \emph{support} of $\varphi$.
Note that the support always contains an even number of vertices.

Using the correspondence between pivot operations
and determinants of submatrices, we can formulate a condition that
specifies the edges present in a graph after a sequence
of pivots.

\begin{Theorem}\label{thm:iterate}
Let $\varphi$  be an applicable sequence of pivoting operations for
$G$, and let $S = \sup(\varphi)$. Then $\det(G\varphi\sub{\pair xy})
= \det(G\sub{S \xor \{x,y\}})$ for $x,y \in V(G)$, $x\neq y$.
Consequently, $\pair xy \in E(G\varphi)$ iff this value equals $1$.
\end{Theorem}
\begin{Proof}
We prove the equality in the statement by induction on the number of
pivot operations.
When $\varphi$ is the empty sequence, we
read the identity
$\det G\sub{\pair xy} = \det G\sub{\varnothing \xor \{x,y\}}$.

So assume $\varphi = [\pair uv] \varphi' $.
Let $S  = \sup(\varphi)$, then $S' = S \xor \{u,v\}$
is the support of $\varphi'$.
We apply the induction hypothesis to
the applicable sequence $\varphi'$ in the
graph $G[\pair uv]$.
Then
$\det G\varphi\sub{\pair xy} =
\det G[\pair uv]\varphi' \sub{\pair xy} =
\det G[\pair uv] \sub{S' \xor \{x,y\}}$.
Now we can apply Theorem~\ref{thm:geelen}, to obtain
$\det G \sub{S' \xor \{x,y\} \xor \{u,v\}}$
which obviously equals
$\det G\sub{S \xor \{x,y\}}$.
\end{Proof}

We now have the following surprising direct consequence of the
previous theorem.

\begin{Theorem}\label{thm:equaldom}
If $\varphi$ and $\varphi'$ are applicable sequences of pivoting
operations for $G$, then $\sup(\varphi) = \sup(\varphi')$ implies
$G\varphi = G\varphi'$.
\end{Theorem}

As a consequence, when calculating the orbit of graphs under the
pivot operation, as done in \cite{danielsen}, we need not consider
every sequence -- only those that have different support.

The next lemma shows, as a direct corollary to
Theorem~\ref{thm:iterate}, that the vertices of the support of an
applicable sequence $\varphi$ induce a subgraph that has a nonzero
determinant.
\begin{Lemma}\label{lem:applic=>odd}
Let $\varphi$ be a sequence of pivots applicable in graph $G$, and
let $S= \sup(\varphi)$.
Then $\det G\sub{S} = 1$.
\end{Lemma}
\begin{Proof}
If $S$ is empty, then indeed $\det G\sub{\varnothing} = 1$.
Now let $S$ (and $\varphi$) be non-empty. Let $\varphi = \varphi'
[uv]$, so $S = \sup \varphi' \xor \{u,v\}$. As $\varphi$ is
applicable, $\pair uv$ must be an edge in $G\varphi'$. By
Theorem~\ref{thm:iterate}, $\det G\varphi'\sub{\pair uv} = \det
G\sub{S} =1$.
\end{Proof}

\bigskip
Two special cases of Theorem~\ref{thm:equaldom} are known from
the literature, the triangle equality (involving three vertices) and
commutativity (involving four vertices).

The \emph{triangle equality} is a classic result in the theory of
pivots. Arratia et al. give a proof \cite[Lemma~10]{arratia}
involving certain graphs with 11 vertices. Independently Genest
obtains this result in his Thesis \cite[Proposition~1.3.5]{genest}.
The cited work of Oum \cite[Proposition~2.5]{oum} contains a proof
which applies Lemma~\ref{lem:oum}.

\begin{Corollary}\label{cor:triangle}
If $u,v,w$ are three distinct vertices in graph $G$
such that $\pair uv$ and $\pair uw$ are edges.
Then $G[uv][vw] = G[uw]$.
\end{Corollary}
\begin{Proof}
Note that
$\pair vw$ is an edge in $G[uv]$ iff
$\det G\sub{ \{v,w\}\xor\{u,v\} } = \det G\sub{ \{w,u\} } = 1$.
The latter holds iff
$\pair uw$ is an edge in $G$.
Hence the pivots at both sides are applicable,
and the result
follows from Theorem~\ref{thm:equaldom}.
\end{Proof}

Another result that fits in our framework is the commutativity of
pivots on disjoint sets of nodes. It was obtained by Harju et al.
\cite{harju:parallelism} (see also \cite{note}) studying graph
operations modelled after gene rearrangements in organisms called
ciliates. The property states that two disjoint pivots $[uv]$ and
$[wz]$, when applicable in either order, have a result independent
of the order in which they are applied.



The next lemma is also proved in Corollary~7 of \cite{nest} using
linear fractional transformations. Essentially, it states that
`twins' stay `twins' after pivoting. Here we obtain it as a
consequence of Theorem~\ref{thm:iterate}.
\begin{Lemma}
Let $v,v'$ be vertices in graph $G$ such that $v \sim_G x$ iff $v'
\sim_G x$ for each vertex $x$. Then for each applicable sequence
$\varphi$ of pivots, $v \sim_{G\varphi} x$ iff $v' \sim_{G\varphi}
x$ for each vertex $x$.
\end{Lemma}
\begin{Proof}
Let $S = \sup(\varphi)$. We have $v \sim_{G\varphi} x$ iff $\det
(G\varphi[\{v\}\xor\{x\}]) = 1$ iff $\det (G[S \xor \{v\}\xor\{x\}])
= 1$ iff $\det (G[S \xor \{v\}\xor\{x\}\xor\{v\}\xor\{v'\}]) = 1$
(since $v \sim_G x$ iff $v' \sim_G x$ for each vertex $x$) iff $\det
(G[S \xor \{v'\}\xor\{x\}]) = 1$ iff $\det
(G\varphi[\{v'\}\xor\{x\}]) = 1$ iff $v' \sim_{G\varphi} x$.
\end{Proof}

\section{Pivots and Perfect Matchings}

There is a direct correspondence between (the parity of) the
determinant of a graph and (the parity of) the number of perfect
matchings in that graph. This correspondence is explained in a paper
by Little \cite{little}, which we essentially follow below. We
include it in our presentation because it allows us to reformulate
some results in terms of a property of the graph itself, rather than
a property of the associated adjacency matrix. We also give an
application, illustrating that the link to perfect matchings adds
some intuition to results from the literature.

We say that a partition $P$ of $V$ is a \emph{pairing of $V$} if it
consists of sets of cardinality two. Let $\tpart(V)$ be the set of
pairings of $V$. A \emph{perfect matching} in $G$ is a pairing $P$
of $V(G)$ such that $P \subseteq E(G)$. Let $\pmatch(G)$ be the
number of perfect matchings of $G$, modulo 2.

For a $V \times V$ matrix $A$, the \emph{Pfaffian} of $A$, denoted
by $\mbox{Pf}(A)$, is defined as $ \sum_{P \in \tpart(V)} \sgn(P)
\prod_{\pair{x}{y} \in P} a_{x,y} $ where $\sgn(P)$ is the sign of a
permutation on the vertices associated with the pairing. As with the
determinants, we apply this notion only for adjacency matrices of
graphs over $\two$, which means $\sgn(P)$ can be dropped from the
formula.

If we evaluate this expression for the adjacency matrix $A$ of a
graph $G$ then we obtain the parity of the number of perfect
matchings the subgraph in $G$: the formula determines, for each
pairing of $V(G)$, whether or not it is a perfect matching.

For skew matrices (where $a_{u,v} = -a_{v,u}$ for all $u,v$)
it is known that $\mbox{Pf}(A)^2 = \det(A)$.
However, over $\two$ every symmetric matrix is skew,
and also the square can be dropped without changing the value.
Thus, for a graph $G$ we know that $\det(G) = \pmatch(G)$.

If we rephrase Theorem~\ref{thm:iterate}
we obtain an elegant characterization of the edges after
pivoting.

\begin{Theorem}\label{thm:again}
Let $\varphi$  be an applicable sequence of pivoting operations for
$G$, and let $S = \sup(\varphi)$.
Then, for vertices $x,y$ with $x\neq y$,
$\pair xy \in E(G\varphi)$ iff
$\pmatch( G\sub{S \xor \{x,y\} } ) = 1$.
\end{Theorem}

For small graphs the number of perfect matchings might be
easier to determine than the determinant.
For instance, for a graph $G$ on four nodes there are
only three pairs of edges that can be present to
contribute to the value $\pmatch(G)$.

\centerline{\begin{picture}(80,20)
\gasset{AHnb=0,Nw=1.5,Nh=1.5,Nframe=n,Nfill=y}
\gasset{ExtNL=y,NLdist=1.5,NLangle=90}
  \node(u1)(05,15){}
  \node(u2)(15,15){}
  \node(v1)(05,05){}
  \node(v2)(15,05){}
  \drawedge(u1,u2){}
  \drawedge(v1,v2){}
  \node(u1)(35,15){}
  \node(u2)(45,15){}
  \node(v1)(35,05){}
  \node(v2)(45,05){}
  \drawedge(u2,v2){}
  \drawedge(v1,u1){}
  \node(u1)(65,15){}
  \node(u2)(75,15){}
  \node(v1)(65,05){}
  \node(v2)(75,05){}
  \drawedge(u2,v1){}
  \drawedge(v2,u1){}
\end{picture}%
}

A commutivity result is obtained in
\cite[Theorem~6.1(iii)]{harju:parallelism}.
Assume $\pair uv$ and $\pair zw$ are edges in $G$ on
four different vertices
$u,v,w,z$. Then both $[uv][wz]$ and $[wz][uv]$ are applicable
iff the induced subgraph $G[\{u,v,w,z\}]$ is not isomorphic to
$C_4$ or $D_4$.

Its proof in \cite{harju:parallelism} is not difficult,
a simple case analysis suffices.
Here we note that $\pair uv$ en $\pair wz$ must be edges
in order for $[uv]$ and $[wz]$ to be applicable.
Both $[uv][wz]$ and $[wz][uv]$ are applicable iff
$\pmatch(G\sub{u,v,w,z})=1$. Thus the subgraph
$G\sub{u,v,w,z}$ must contain either one or three perfect matchings,
where the first $\{\pair uv, \pair wz\}$ is given.
Two perfect matchings occur precisely when the subgraph is isomorphic
to $C_4$ or $D_4$.

\centerline{\begin{picture}(50,20)
\gasset{AHnb=0,Nw=1.5,Nh=1.5,Nframe=n,Nfill=y}
\gasset{ExtNL=y,NLdist=1.5,NLangle=90}
  \node(u1)(05,15){}
  \node(u2)(15,15){}
  \node(v1)(05,05){}
  \node(v2)(15,05){}
  \drawedge(u1,u2){}
  \drawedge(u2,v2){}
  \drawedge(v2,v1){}
  \drawedge(v1,u1){}
  \put(18,05){$C_4$}
  \node(u1)(35,15){}
  \node(u2)(45,15){}
  \node(v1)(35,05){}
  \node(v2)(45,05){}
  \drawedge(u1,u2){}
  \drawedge(u2,v2){}
  \drawedge(v2,v1){}
  \drawedge(v1,u1){}
  \drawedge(v1,u2){}
  \put(48,05){$D_4$}
\end{picture}%
}

As was noted below Theorem~\ref{thm:geelen}, we may also look for a
non-empty set $S$ such that every $v \in V(G)$ is adjacent to an
even number of vertices of $S$. E.g., for $D_4$ we can take $S$ to
be the set of the two vertices that are not connected by an edge.

\section{Reduced Sequences}

We have seen that if we have an applicable sequence of pivots, then
the result of that series of operations only depends on the support,
the set of vertices occurring an odd number of times as a
pivot-vertex. This does not automatically mean that the sequence can
be reduced to an equivalent sequence in which each vertex occurs
only once. This because one needs to verify that all the operations
are applicable, i.e., that all pivot-pairs are edges in the graph to
which they are applied.

We call a sequence of pivots \emph{reduced} \cite{genest} if no
vertex occurs more than once in the pivots. It turns out that we can
use a greedy strategy to reduce a sequence of pivot operations,
growing a sequence with given support.

Let $G$ be a graph, and let $S\subseteq V(G)$ be the support
of an applicable sequence of
pivots for $G$.
We will construct a reduced sequence of pivots
with support $S$.
Obviously we may assume that $S$ is non-empty.
Observe that since $\det G\sub {S} =1$,
there must be at least one element in the adjacency matrix
of $G$ that is non-zero, i.e., there is an edge $\pair uv$
in $G\sub {S}$.

Apply $[uv]$ to graph $G$, and proceed iteratively with graph
$G[uv]$ and support set $S-\{u,v\}$, where $\det G[uv]
\sub{D-\{u,v\}} = 1$ again holds (by Theorem~\ref{thm:geelen}) and
we stop when we have exhausted the support.

\begin{Theorem}
For every applicable sequence of pivots $\varphi$
there exists an
applicable reduced sequence $\varphi'$
such that $\sup(\varphi) = \sup(\varphi')$
--- and therefore $G\varphi = G\varphi'$.
\end{Theorem}

\begin{Remark}
The possibility to construct an applicable reduced sequence
with given support depends on the fact that there must
be at least one edge to obtain a non-zero determinant.
In fact every column in the matrix must contain at least
one edge. This means we can even choose one of the
vertices of the pivot.

As an example, we return to the topic of commutivity.
It is known that if $[uv][wz]$ is
applicable, then we can not conclude that $[wz][uv]$ is applicable.
However, $\det G\sub {u,v,w,z} =1$, so
we can construct an applicable sequence with
support $\{u,v,w,z\}$.
Fixing $z$ we know that there is an edge adjacent to that vertex,
which can be either $\pair wz$, $\pair vz$ or $\pair uz$.
When pivoting over this edge, the remaining two vertices
must form an edge in the graph.

Hence, we have shown the following fact:
if, for for different vertices, $[uv][wz]$
is applicable, then either at least one of the pivot sequences
$[wz][uv]$, $[vz][uw]$, or $[uz][vw]$ is applicable.
This is essentially Lemma~1.2.11 of \cite{genest}.
\end{Remark}

The previous theorem shows that also the converse of
Lemma~\ref{lem:applic=>odd} holds.

\begin{Theorem}\label{thm:applic<=>odd}
Let $S$ be a set of vertices of graph $G$. Then $\det G\sub S = 1$
iff there exists a (reduced) sequence of pivots $\varphi$  with
support $S$ that is applicable in $G$.
\end{Theorem}

The size of $\{ S \subseteq V \mid \det G\sub S = 1 \}$ is precisely
the value of the interlace polynomial $q(G)$ of $G$ on $x=1$, see
\cite[Corollary~2]{aigner}.

\section{Graphs with Self-Loops} \label{sec_self_loops}
Until now we have considered simple graphs (graphs without loops or
parallel edges). In this section we consider graphs $G$ with loops
but without parallel edges. The adjacency matrices $A$ corresponding
to such graphs are precisely the symmetrical $(0,1)$-matrices. If
vertex $u$ has a loop in $G$, then the matrix $A \sub{\{u\}}$ is
equal to the $1 \times 1$ matrix $(1)$. Hence, $\det A*\{u\} = 1$
and the general pivot of Section~\ref{sec_pivots} is defined, and is
modulo $2$ equal to
$$
A*\{u\} = \left(
\begin{array}{c|c}
1 & \chi_u^T \\
\hline \chi_u & A\sub{V-u} - \chi_u \chi_u^T
\end{array}
\right),
$$
where $\chi_u$ is the column vector belonging to $u$ without element
$a_{uu}$. We define the elementary pivot $G*u$ for loop vertex $u$
on $G$ by the graph corresponding to adjacency matrix $A*\{u\}$. The
elementary pivot $G*u$ is obtained from $G$ by complementing the
neighbourhood $N_G(u)$ of $u$ (just as in simple graphs) \emph{and},
for $v \in N_G(u)$, we add a loop to $v$ if $v$ is a non-loop vertex
in $G$ and remove the loop if $v$ is a loop vertex in $G$. Hence, we
will call $G*u$ \emph{local complementation} (on graph $G$ with loop
$u$). We can apply Proposition~\ref{prop:geelen}, and similar to
Theorem~\ref{thm:geelen} we obtain (in $\two$) the following result.

\begin{Theorem}\label{thm:geelen2}
Let $G$ be a graph, and let $u \in V(G)$ be a vertex that has a
loop. Then, for $Y \subseteq V(G)$,
\[
    \det ((G*u)\sub{Y}) = \det (G\sub{Y \xor \{u\}})
\]\qed
\end{Theorem}

The pivot operation on edge $e$ for graph with loops is identical to
that operation for simple graphs: it is only defined if both
vertices of $e$ do not have loops and it does not remove or add any
loop of the graph.

Results of the previous sections carry over to sequences $\varphi$
of operation having both local complementation and pivot operations.
In particular, Theorems~\ref{thm:equaldom} and
~\ref{thm:applic<=>odd} carries over.

\begin{Theorem}\label{thm:equaldom2}
If $\varphi$ and $\varphi'$ are applicable sequences for $G$ having
(possibly) both local complementation and pivot operations, then
$\sup(\varphi) = \sup(\varphi')$ implies $G\varphi = G\varphi'$.
Also, let $S \subseteq V(G)$. Then $\det G\sub S = 1$ iff there
exists a (reduced) sequence $\varphi$ with support $S$, having
(possibly) both local complementation and pivot operations, that is
applicable in $G$.
\end{Theorem}

The size of $\{ S \subseteq V \mid \det G\sub S = 1 \}$, for graph
$G$ with loops, is precisely the value of a polynomial $Q(G)$,
defined in \cite[Section~4]{aigner} and related to the interlace
polynomial $q(G)$, of $G$ on $x=2$, see \cite[Corollary~5]{aigner}.
Moreover, the previous theorem can also be stated in terms of
\emph{general perfect matchings}: considering a loop on $x$ as the
edge $\{x\} \in E(G)$, then a general perfect matching is a $P
\subseteq E(G)$ that is a partition of $V(G)$.

\begin{Remark}
In the theory of gene assembly in ciliates\cite{ciliates}, the local
complementation operation on $u$ with the removal of $u$ is called
\emph{graph positive rule}, and the pivot operation on $uv$ with the
removal of both $u$ and $v$ is called \emph{graph double rule}.
These rules are defined on \emph{signed graphs}, where each vertex
is labelled by either $-$ or $+$. Now, label $-$ corresponds to a
non-loop vertex and label $+$ corresponds to a loop vertex. Hence,
we obtain the result that any two sequences of these graph rules
with equal support obtain the same graph. Moreover, we obtain that a
signed graph can be transformed into the empty graph by these graph
rules iff the determinant of corresponding adjacency matrix has
determinant $1$ modulo $2$.
\end{Remark}

\section{Discussion}
We have related applicable sequences of pivot operations to
determinants and perfect matchings in a graph. In this way, we have
shown that two applicable sequences of pivot operations with equal
support have the same effect on the graph. Moreover, for a given set
$S$ of vertices, we have shown that there is a applicable sequence
$\varphi$ of pivot operations with support $S$ precisely when the
number of perfect matchings of the subgraph induced by $S$ is odd
(or equivalently, when the determinant of the adjacency matrix of
the subgraph is odd). In fact, there is an applicable reduced
sequence $\varphi'$ with equal support as $\varphi$. Finally, we
have shown that pivots and local complementation can `work together'
in the case of graphs with loops in the sense that equal support
renders equal graphs.


\newpage
\appendix

\section{Pivots and Matchings}

In this appendix we give an independent proof of
Theorem~\ref{thm:again} in the style of Oum \cite{oum}, using
perfect matchings instead of determinants, as it may be of
independent interest. The proof was made superfluous when the
authors discovered references \cite{little} and \cite{geelen}.
However, the proofs in this appendix are straightforward and
therefore the reader may prefer this approach.

\bigskip
Recall that $x \sim_G y$ if either $\pair xy \in E(G)$ or $x=y$.

As a technical tool we need a formula that can be used to compute
the number of perfect matchings in a graph, but which can also be
applied when we have duplicate vertices.

For (an even number of) variables $x_1,\dots, x_n$ let
$\pmatch_G(x_1,\dots,x_n)$ denote the following
logical expression:
\[
\bigoplus_{P \in \tpart \{x_1,\dots,x_n\} }
\bigwedge_{\pair{x}{y} \in P} (x \sim_G y)
\]
The number of variables used in the expression varies; we
assume this number is clear from the context. Clearly,
$\pmatch_G(x,y)$ equals $x \sim_G y$. Moreover $\pmatch_G()$ is true
-- the logical and $\wedge$ over $0$ arguments is (considered) true,
and the logical exclusive or over $1$ argument $a$ is (considered)
$a$. This is in line with the fact that there is a single perfect
matching on zero vertices.

If we evaluate this expression for the (pairwise different) vertices
$v_1,\dots,v_n$ of graph $G$ then we obtain the value
$\pmatch_G(v_1,\dots,v_n)$ which equals
$\pmatch(G\sub{\{v_1,\dots,v_n\}})$, the parity of the number of
perfect matchings the subgraph in $G$ induced by $v_1,\dots,v_n$
(identifying 0 and 1 with false and true, respectively). Due to the
highly symmetric form of the formula $\pmatch_G$ the ordering of the
vertices as arguments to the formula is not important for the value.
We will use this fact frequently below.

The formula can also be evaluated when two (or more) of its
arguments are chosen to be the same vertex in the graph. The next
result shows equal vertices can be omitted (in pairs).

\begin{Lemma}\label{lem:equivmat}
Let $v_1,\dots,v_{n-2},v,v'$ be vertices in graph $G$ such that
$v \sim_G x$ iff $v' \sim_G x$ for each vertex $x$.
Then $\pmatch_G(v_1,\dots,v_{n-2},v,v') = \pmatch_G(v_1,\dots,v_{n-2})$.
\end{Lemma}
\begin{Proof}
Observe that the condition of the lemma on $v$ and $v'$
implies that $v\sim_G v'$ holds.

For $n=2$ the left hand side $\pmatch_G(v,v') $
equals $v\sim_G v'$ which equals the right hand side
$\pmatch_G() $ which has been set to true.

Now let $n>2$. In the formal expression $\pmatch_G$ each pairing $P$
that does not contain $\pair{x_{n-1}}{x_n}$ has two pairs $\pair{x_{n-1}}{x_i}$ and
$\pair{x_n}{x_j}$. For $P$ there is a (unique) $P'$ corresponding to $P$
where $P' \backslash P = \{\pair{x_{n-1}}{x_j}, \pair{x_n}{x_i} \}$ (and hence $P
\backslash P' = \{\pair {x_{n-1}}{x_i}, \pair{x_n}{x_j} \}$).
Since $v \sim_G x$ iff
$v' \sim_G x$ for each vertex $x$, we have
$(v \sim_G x_i) \land (v' \sim_G x_j) = (v \sim_G x_j) \land (v' \sim_G x_i)$.
Hence the contributions of pairings $P$ and $P'$ cancel.

The remaining pairings all contain $\pair{v}{v'}$ for which $v
\sim_G v'$ can be dropped from the formula, as $v \sim_G v'$ holds.
The resulting formula equals that of $\pmatch_G(v_1,\dots,v_{n-2})$.
\end{Proof}

The next lemma shows that we can characterize pivoting by the parity
of the number of perfect matchings in subgraphs.
It is a simple reformulation of the result of Oum,
but essential as a first step to understand
the connection between pivoting and perfect matchings.

\begin{Lemma}\label{lem:oumagain}
Let $G = (V,E)$ be a graph, and fix $\pair uv \in E$. For $x,y \in
V$ we have
$\pmatch_{G[uv]}(x,y) = \pmatch_G(x,y,u,v)$.
\end{Lemma}
\begin{Proof}
In the evaluation of $\pmatch_{G[uv]}(x,y)$
we consider a single pair $x \sim_{G[uv]} y$ only,
the left-hand side in Lemma~\ref{lem:oum}.

As $u \sim_G v$ holds, we may replace the factor $x \sim_G y$
in the statement of Lemma~\ref{lem:oum} by $(x \sim_G y \land u \sim_G v)$.
Now the right-hand side of the formula equals
$\pmatch_G(x,y,u,v)$.
\end{Proof}

The main technical result is a generalization of the
previous lemma, which now includes an additional
sequence of nodes on both sides.
Before stating this result we explicitly compute the simplest
of these generalizations, with variables $x_1,x_2,x_3,x_4$
instead of $x,y$.
This example visualizes the more general argumentations in
the proof of our general result, which follows the example.


\begin{Example}
$\pmatch_{G[uv]}(x_1,x_2,x_3,x_4)$
equals
$
 (x_1 \sim_{G[uv]} x_2 \land x_3 \sim_{G[uv]} x_4)  \xor
 (x_1 \sim_{G[uv]} x_3 \land x_2 \sim_{G[uv]} x_4)  \xor
 (x_1 \sim_{G[uv]} x_4 \land x_2 \sim_{G[uv]} x_3)
$.

Now substitute each $x\sim_{G[uv]}y$ by the formula given in
Lemma~\ref{lem:oum},
to obtain

\noindent
$
  (\; [x_1x_2 \xor (x_1u\land x_2v) \xor (x_1v\land x_2u) ]   \land
        [x_3x_4 \xor (x_3u\land x_4v) \xor (x_3v\land x_4u) ]
  \;)
\bigoplus
  (\; [x_1x_3 \xor (x_1u\land x_3v) \xor (x_1v\land x_3u) ]   \land
        [x_2x_4 \xor (x_2u\land x_4v) \xor (x_2v\land x_4u) ]
  \;)
\bigoplus
  (\; [x_1x_4 \xor (x_1u\land x_4v) \xor (x_1v\land x_4u) ]   \land
        [x_2x_3 \xor (x_2u\land x_3v) \xor (x_2v\land x_3u) ]
  \;)
$,
\\
where we write $xy$ rather than $x\sim_G y$.

By distributivity (i.e. using the logical identity $a \wedge (b \xor
c) = (a \wedge b) \xor (a \wedge c)$) this is equivalent to

\noindent
$
   ( x_1x_2 \land x_3x_4 ) \xor
   ( x_1x_2 \land x_3u \land x_4v ) \xor
   ( x_1x_2 \land x_3v \land x_4u ) \xor
   ( x_1u \land x_2v \land x_3x_4 ) \xor
   ( x_1u \land x_2v \land x_3u \land x_4v ) \xor
   ( x_1u \land x_2v \land x_3v \land x_4u ) \xor
   ( x_1v \land x_2u \land x_3x_4 ) \xor
   ( x_1v \land x_2u \land x_3u \land x_4v ) \xor
   ( x_1v \land x_2u \land x_3v \land x_4u )
\bigoplus
   ( x_1x_3 \land x_2x_4 ) \xor
   ( x_1x_3 \land x_2u \land x_4v ) \xor
   ( x_1x_3 \land x_2v \land x_4u ) \xor
   ( x_1u \land x_3v \land x_2x_4 ) \xor
   ( x_1u \land x_3v \land x_2u \land x_4v ) \xor
   ( x_1u \land x_3v \land x_2v \land x_4u ) \xor
   ( x_1v \land x_3u \land x_2x_4 ) \xor
   ( x_1v \land x_3u \land x_2u \land x_4v ) \xor
   ( x_1v \land x_3u \land x_2v \land x_4u )
\bigoplus
   ( x_1x_4 \land x_2x_3 ) \xor
   ( x_1x_4 \land x_2u \land x_3v ) \xor
   ( x_1x_4 \land x_2v \land x_3u ) \xor
   ( x_1u \land x_4v \land x_2x_3 ) \xor
   ( x_1u \land x_4v \land x_2u \land x_3v ) \xor
   ( x_1u \land x_4v \land x_2v \land x_3u ) \xor
   ( x_1v \land x_4u \land x_2x_3 ) \xor
   ( x_1v \land x_4u \land x_2u \land x_3v ) \xor
   ( x_1v \land x_4u \land x_2v \land x_3u )
$

There are twelve terms with four factors, which are six different
terms each occurring twice, hence cancelling each other. The three
terms with two factors can be extended adding a third term $uv$
(which is true). Rearranging these 15 remaining terms we get

\noindent
$
   ( x_1x_2 \land x_3x_4 \land uv ) \xor
   ( x_1x_2 \land x_3u \land x_4v ) \xor
   ( x_1x_2 \land x_3v \land x_4u ) \xor
   ( x_1x_3 \land x_2x_4 \land uv ) \xor
   ( x_1x_3 \land x_2u \land x_4v ) \xor
   ( x_1x_3 \land x_2v \land x_4u ) \xor
   ( x_1x_4 \land x_2x_3 \land uv ) \xor
   ( x_1x_4 \land x_2u \land x_3v ) \xor
   ( x_1x_4 \land x_2v \land x_3u ) \xor
   ( x_1u \land x_2v \land x_3x_4 ) \xor
   ( x_1u \land x_3v \land x_2x_4 ) \xor
   ( x_1u \land x_4v \land x_2x_3 ) \xor
   ( x_1v \land x_2u \land x_3x_4 ) \xor
   ( x_1v \land x_3u \land x_2x_4 ) \xor
   ( x_1v \land x_4u \land x_2x_3 )
$

These happen to be the fifteen pairings
making up
$\pmatch_G(x_1,x_2,x_3,x_4,u,v)$.
\end{Example}

As announced, the proof of our general result follows the path
sketched in the previous example. It is the `perfect matching
counterpart' of Theorem~\ref{thm:geelen}.

\begin{Theorem}\label{thm:basic}
Let $G$ be a graph, let $v_1,\dots,v_n$ be vertices in $G$, let
$\pair uv \in E(G)$. Then $\pmatch_{G[uv]}(v_1,\dots,v_n) =
\pmatch_G(v_1,\dots,v_n,u,v)$.
\end{Theorem}
\begin{Proof}
If $n=0$ the left hand side equals $\pmatch_{G[uv]}()$
which is true, while the right hand side
$\pmatch_G(u,v)$ is equivalent to $u \sim_G v$,
which is also true, as $\pair uv$ is an edge in $G$.

Now let $n\ge 2$.
For $\pmatch_{G[uv]}(v_1,\dots,v_n)$ the following formula has to be evaluated
\[
\bigoplus_{P \in \tpart\{x_1,\dots,x_n\}}
\bigwedge_{\pair{x}{y} \in P} (x \sim_{G[uv]} y)
\]
According to Lemma~\ref{lem:oum} the relation $\sim_{G[uv]}$ can be replaced
by a suitable expression involving $\sim_G$ in the original graph $G$.
\[
\bigoplus_{P \in \tpart\{x_1,\dots,x_n\}}
\bigwedge_{\pair xy \in P}
 \left( (x \sim_G y) \xor
 (x \sim_G u \land y \sim_G v) \xor
 (x \sim_G v \land y \sim_G u) \right)
\]
Now, we apply the logical identity
$a \wedge (b \xor c) = (a \wedge b) \xor (a \wedge c)$
iteratively to the inner part $\bigwedge_{\pair xy \in P} (\dots)$,
and we obtain for each
$P \in \tpart\{x_1,\dots,x_n\}$ the exclusive or over
a total of $3^{n/2}$ terms,
each of which is a conjunction of factors
of one of the forms
 $x \sim_G y$,
 $(x \sim_G u \land y \sim_G v)$ and
 $(x \sim_G v \land y \sim_G u)$.
Moreover in each such term the variables
$x_1,\dots,x_n$ each occur exactly once.

Now consider such a term in which the constant $u$
occurs $k$ times paired to $x_{i_1}, \dots, x_{i_k}$,
which implies also $v$ occurs $k$ times
paired to certain $x_{j_1}, \dots, x_{j_k}$.

Up to the order of factors,
this term is present in the list
that belongs to any $P'$ that pairs the variables
$x_{i_1}, \dots, x_{i_k}$,
to the $x_{j_1}, \dots, x_{j_k}$
(in any combination)
and equals $P$ for the other variables.
There are $k!$ such pairings, thus $k!$ copies of equivalent terms.
These copies cancel if $k!$ is even, which means if $k\ge 2$.

Hence, for each $P$ we need only consider those terms
for which there is at
most one occurrence of both $u$ and $v$.
Thus we have reduced the previous equation to
\[\begin{array}{cl}
\bigoplus_{P \in \tpart\{x_1,\dots,x_n\}}
&
 \left( \bigwedge_{\pair{x_1}{y_1} \in P}
x_1 \sim_G y_1 \right) \oplus
\\
&
\bigoplus_{\pair{x_1}{y_1} \in P}
\left[ \left( x_1 \sim_G u \wedge
y_1 \sim_G v \wedge
\bigwedge_{\pair xy \in P \backslash \pair{x_1}{y_1}}
x \sim_G y \right)
\right.
\\
&
\left.
\oplus \left( x_2 \sim_G u \wedge x_1 \sim_G v \wedge
\bigwedge_{\pair xy \in P \backslash \pair{x_1}{x_2}} x \sim_G y
\right) \right]
\end{array}
\]

Because $\bigwedge_{\pair xy \in P} x \sim_G y =
\bigwedge_{\pair xy \in P \cup \{\pair uv\}} (x \sim_G y)$,
this is equivalent to

$$\bigoplus_{P \in \tpart \{x_1,\dots,x_n,u,v\} } \left(
\bigwedge_{\pair xy \in P} x \sim_G y \right)$$
and this in turn is the expression that has to be evaluated
for $\pmatch_G(x_1,\dots,x_n,u,v)$.
\end{Proof}

By Lemma~\ref{lem:equivmat}, the previous theorem may be rephrased
as follows, cf. Theorem~\ref{thm:geelen}.
\begin{Theorem}\label{thm:basic2}
Let $G$ be a graph, and let $\pair uv \in E(G)$. Then, for $Y
\subseteq V(G)$,
\[
    \pmatch ((G[\pair uv])\sub{Y}) = \pmatch (G\sub{Y \xor \{u,v\}})
\]\qed
\end{Theorem}
The results of Section~\ref{sec_seq_pivots} involving $\det(G)$ can
hence also be developed using $\pmatch(G)$ through
Theorem~\ref{thm:basic2}. Hence, we obtain, e.g.,
Theorem~\ref{thm:again}.

The following special case of our general result
Theorem~\ref{thm:basic} is a reformulation in the style of the
original Lemma~\ref{lem:oum}, summing over edges in the subgraph of
$G$ induced by $\{u,v,w,z\}$ (with some care in the case of multiple
occurrences of vertices).

\begin{Theorem}
If $[uv][wz]$ is applicable to $G$, then
\[
x \sim_{G[uv][wz]} y = x\sim_G y \bigoplus_{\mbox{\shortstack{
$\{\pair{x_1}{x_2},\pair{x_3}{x_4} \}$\\
 pairing of\\  $\{u,v,w,z\}$ \\ with $ x_1\sim_G x_2$ }}}
((x \sim_G x_3) \wedge (y \sim_G
x_4)) \xor ((x \sim_G x_4) \wedge (y \sim_G x_3))
\]
\end{Theorem}
\begin{Proof}
The result is obtained by rewriting the expression for
$\pmatch_G(x,y,u,v,w,z)$, and using the fact that
$\pmatch_G(u,v,w,z)$ holds.
\end{Proof}

\end{document}